\documentclass[a4paper]{article}

\usepackage[english]{babel}
\usepackage[utf8x]{inputenc}
\usepackage[T1]{fontenc}

\usepackage[top=3cm,bottom=2cm,left=3cm,right=3cm,marginparwidth=1.75cm]{geometry}
\usepackage{authblk}
\usepackage{blindtext}
\usepackage{amsmath, amsthm, amssymb, esint, mathtools, mathabx}
\usepackage{amsfonts}
\usepackage{graphicx}
\graphicspath{{/texfiles/}}
\usepackage[colorinlistoftodos]{todonotes}
\usepackage[colorlinks=true, allcolors=blue]{hyperref}
\usepackage{multicol}

\numberwithin{equation}{section}

\title{$L^2\times L^2\times L^2\to L^{2/3}$ boundedness for trilinear multiplier operator}

\author{A. Martina Neuman}
\affil{Department of Mathematics, New York University, Shanghai}
\affil[ ]{\textit{marsneuman@nyu.edu}}

\begin{document}

\maketitle

\begin{abstract} 

\noindent
This paper discusses the boundedness of the trilinear multiplier operator $T_{m}(f_1,f_2,f_3)$, when the multiplier satisfies a certain degree of smoothness but with no decaying condition and is $L^{q}$-integrable with an admissible range of $q$. The boundedness is stated in the terms of $\|m\|_{L^{q}}$. In particular, 
\begin{equation*}\|T_{m}\|_{L^2\times L^2\times L^2\to L^{2/3}}\lesssim\|m\|_{L^{q}}^{q/3}.\end{equation*}
\end{abstract}

%%  LaTeX will not make the title for the paper unless told to do so.

\tableofcontents

\section{Introduction}

\noindent
Let $d\geq 1$. Let $m(\xi,\eta,\delta)$ be a function on $\mathbb{R}^{3d}$. Define an operator $T_{m}$ as follows.
\begin{equation*}T_{m}(f,g,h)(x) = \int_{\mathbb{R}^{d}}\int_{\mathbb{R}^{d}}\int_{\mathbb{R}^{d}} m(\xi,\eta,\delta)\hat{f}(\xi)\hat{g}(\eta)\hat{h}(\delta)e^{2\pi ix\cdot(\xi+\eta+\delta)}\,d\xi d\eta d\delta.\end{equation*}
This paper aims to give an explicit result on the $L^2\times L^2\times L^2\to L^{2/3}$ boundedness of $T_{m}$ given (restricted) smoothness and integrability of $m$. The author anticipates the use of such result in an upcoming project. At the time of starting this paper, the author was not aware of any such explicit result on trilinear multiplier operator on the market. The crucial point here, that will be needed for later use, is the dependence of either operator norm $\|T_{m}\|$ on $\|m\|_{L^{q}}$ (see $\ref{eq 5}$). Unfortunately, due to the lack of duality theory for $L^{s}, 0<s<1$, one can't extend easily this result to other exponents in the Banach range, $L^{p_1}\times L^{p_2}\times L^{p_3}\to L^{r}$, with $1/p_1+1/p_2+1/p_3 = 1/r$. Such obstacle is not met in the case of bilinear setting; see [7].\\

\noindent
There is a body of literature regarding the boundedness of multiplier operator in the linear and bilinear settings, with various conditions on $m$, ranging from decay to smoothness to integrability. In the bilinear setting, one classical condition on $m$ to guarantee the boundedness of $T_{m}$ on the Banach range is the Coifman-Meyer condition [2]:
\begin{equation*}|\partial^{\alpha}m(\xi,\eta)|\leq C_{\alpha}|(\xi,\eta)|^{-|\alpha|}\end{equation*}
for sufficiently many $\alpha$. Moreover, in the bilinear setting, if one is to merely impose uniform bounds on derivatives of $m$: $\|\partial^{\alpha}m\|_{L^{\infty}}\leq C_0<\infty$, then one needs to make other compromises. It was shown in [1] that such uniform boundedness on $m$ alone is not sufficient to make $T_{m}$ a bounded operator on $L^2\times L^2\to L^1$. It was shown in [6] that if one further imposes $L^2$-integrability of $m$, one can get back a bounded operator. Moreover, the same authors in [6] showed that there is a short range of integrability that one can impose in order to secure boundedness of $T_{m}$. This paper follows the ideas in [6]. That means, the multipliers considered here only have uniform derivative bounds plus some integrability.\\

\noindent
There are other venues in the multilinear settings, where positive boundednesss results have been established when the multiplier smoothness is much compromised. See [8] for the case of multipliers from Sobolev spaces, and [5] for the case of $L^{r}$-based Sobolev spaces.

\section{Notation explanation}

\noindent
$\mathbb{N}_0$: the set of nonnegative integers\\

\noindent
$\ominus$: orthogonal complement\\

\noindent
$\{F,M\}^{d*}$: $d$-tuples whose elements are either $F$ or $M$ and at least one of those must be $M$\\

\noindent
$\mathcal{C}^{J}_{c}$: the space of all continuous functions whose derivatives up to, and including, order $J$th are continuous and which have compact support\\

\noindent
$\kappa Q$: a cube (interval) that has the same center as the cube (interval) $Q$ and $\kappa$ times the side-length of $Q$\\

\noindent
$|\cdot |$: either means an absolute value or the Euclidean norm of a vector or the cardinality of a discrete set

\section{Some multiresolution analysis background}

\noindent
The approach followed in this discussion requires an understanding of multiresolution analysis and wavelets. In this section, necessary background is introduced. The following facts can be found in [10].\\ 

\noindent
First, one starts with an orthonormal basis of $L^2(\mathbb{R})$.\\

\noindent
{\bf Definition 1.} An (inhomogeneous) multiresolution analysis is a sequence $\{V_{j}: j\in\mathbb{N}_0\}$ of subspaces of $L^2(\mathbb{R})$ such that\\
a) $V_0\subset V_1\subset\cdots\subset V_{j}\subset V_{j+1}\subset\cdots$ spans $\bigcup_{j\geq 0}V_{j}=L^2(\mathbb{R})$.\\
b) $f\in V_0$ iff $f(x-n)\in V_0$ for any $n\in\mathbb{Z}$.\\
c) $f\in V_{j}$ iff $f(2^{-j}\cdot)\in V_0$ for $j\in\mathbb{N}$.\\
d) There exists $\phi_{F}\in V_0$ such that $\{\phi(\cdot -n)\}_{n\in\mathbb{Z}}$ is an orthonormal basis in $V_0$.\\

\noindent
Because of the properties $(c), (d)$, $\phi_{F}$ is called the {\it scaling function} or the {\it father wavelet} of such system.\\

\noindent
For $j\in\mathbb{N}_0$, let $W_{j}=V_{j+1}\ominus V_{j}$. Let $\phi_{M}\in W_0$ be such that $\{\phi_{M}(\cdot - n)\}_{n\in\mathbb{Z}}$ forms an orthonormal basis in $W_0$. Such existence is guaranteed by wavelet theory [10]. This function $\phi_{M}$ is called the {\it mother wavelet} associated with $\phi_{F}$.\\

\noindent
One now needs an orthonormal system of $L^2(\mathbb{R}^{d})$. Such system can be generated from one of one dimension.\\

\noindent
For $n = (n_{r})_{1\leq r\leq d}\in\mathbb{Z}^{d}$, let 
\begin{equation*}\Phi_{n}(x) = \prod_{r=1}^{d}\phi_{F}(x_{r}-n_{r}),\end{equation*} 
with $x\in\mathbb{R}^{d}$.\\

\noindent
Denote $G=(G_1,\cdots,G_{d})\in\{F,M\}^{d}$ - in other words, a $d$-tuple of types, $F$ or $M$ (father wavelet or mother wavelet, respectively). Let
\begin{equation*}\Phi^{G}_{n}(x) = \prod_{r=1}^{d}\phi_{G_{r}}(x_{r}-n_{r}).\end{equation*} 

\noindent
Finally, let $G^0=\{(F,\cdots,F)\}$ and $G^{j}=\{F,M\}^{d*}$ for all $j\in\mathbb{N}$. Assume that $\|\phi_{F}\|_{L^2}=\|\phi_{M}\|_{L^2}=1$, then one has the following promise.\\

\noindent
{\bf Proposition 2.} The system below forms an orthonormal basis in $L^2(\mathbb{R}^{d})$,
\begin{equation}\label{eq 1}\Phi_{n}^{j,G}(x) = \begin{cases}\Phi_{n}(x) &\text{ for } j=0, G\in G^0, n\in\mathbb{Z}^{d}\\ 2^{\frac{jd}{2}}\Phi_{n}^{G}(2^{j}x) &\text{ for } j\in\mathbb{N}, G\in G^{j}, n\in\mathbb{Z}^{d} \end{cases}.\end{equation}

\noindent
{\it Remark 3:} One can also require that $\phi_{F},\phi_{M}\in\mathcal{C}^{J}_{c}$ and that $\int x^{\alpha}\phi_{M}(x)\,dx=0$ for all $|\alpha|\leq K$; here $J, K$ are sufficiently large positive integers for the following calculations to hold. For the guarantee of smoothness and moment cancellations, see Theorem 1.61 and remark 1.62 in [9].\\

\noindent
In what follows one is concerned with functions on $\mathbb{R}^{3d}$ is in the place of $d$ - hence for example, $d$ is replaced by $3d$ in $\ref{eq 1}$. \\

\noindent
{\it Remark 4:} By definition of $\Phi^{j,G}_{n}$ in $\ref{eq 1}$, each $\Phi^{j,G}_{n}$ with $j\in\mathbb{N}_0, n\in\mathbb{Z}^{3d}$ can be written as $\Phi^{j,G}_{(k_1,k_2,k_3)}$ with $k_{i}\in\mathbb{Z}^{d}$, and moreover, $\Phi^{j,G}_{(k_1,k_2,k_3)}=\omega_{1,k_1}\cdots\omega_{l,k_{l}}$ with $\omega_{i,k_{i}}$'s being functions of only variables $(x_{(i-1)d+1},\cdots,x_{id})$. By the assumption made in {\it Remark 3}, $\omega_{1,k_1},\omega_{2,k_2},\omega_{3,k_3}\in\mathcal{C}^{J}_{c}$, for some sufficiently large $J$. For a fixed $j,G$, the supports of $\Phi^{j,G}_{(k_1,k_2,k_3)}$ are disjoint. Hence in particular, he supports of $\omega_{i,k_{i}}$ have finite overlaps in $k_{i}$. Moreover 
\begin{equation}\label{eq 2}
\|\omega_{1,k_1}\|_{L^{\infty}}, \|\omega_{2,k_2}\|_{L^{\infty}},\|\omega_{3,k_3}\|_{L^{\infty}}\lesssim_{d} 2^{\frac{jd}{2}}.\end{equation}
An example of how the facts above are (frequently) used in the following analysis, is as follows,
\begin{equation*}
\sum_{k_1}\|\omega_{1,k_1}\hat{f}_1\|_{L^2}^2\lesssim_{d}\|\sum_{k_1}\omega_{1,k_1}\hat{f}_1\|_{L^2}^2\leq\|f_1\|_{L^2}\|\sum_{k_1}\omega_{1,k_1}\|_{L^{\infty}}\lesssim_{d} 2^{jd/2}\|f_1\|_{L^2}.\end{equation*}
Here the sufficient disjointness in supports of $\omega_{1,k_1}$'s is used in the first inequality, and $\ref{eq 2}$ is used in the last.\\

\noindent
Let $m(\xi_1,\cdots,\xi_{l})$ be a function on $\mathbb{R}^{3d}$ where $\xi_1,\cdots,\xi_{l}\in\mathbb{R}^{d}$. The following lemma is essentially given in [6].\\

\noindent
{\bf Lemma 3.} Let $K$ be a positive integer. Assume that $m\in\mathcal{C}^{K+1}$ is a function on $\mathbb{R}^{3d}$ such that
\begin{equation}\label{eq 3} \sup_{|\alpha|\leq K+1}\|\partial^{\alpha}m\|_{L^{\infty}}\leq C_0<\infty.\end{equation}
Then one has,
\begin{equation}\label{eq 4} |\langle\Phi^{j,G}_{n},m\rangle|\lesssim_{C_0} 2^{-(K+1+d)j}\end{equation}
provided that $\phi_{M}$ has $K$ vanishing moments.\\

\noindent
{\it Remark 5:} See the appendix for a brief discussion about the proof of this lemma.

\subsection{Main theorem}

\noindent
{\bf Theorem 4.} Let $1\leq q<3$. Set $M_{q}=\lfloor\frac{3d}{3-q}\rfloor +1$. Let $m(\xi_1,\xi_2,\xi_3)$ be a function in $L^{q}(\mathbb{R}^{3d})\cap\mathcal{C}^{M_{q}}(\mathbb{R}^{3d})$ satisfying
\begin{equation*} \|\partial^{\alpha}m\|_{L^{\infty}}\leq C_0 \end{equation*}
for all $|\alpha|\leq M_{q}$. Then the trilinear operator $T_{m}$ with the multiliplier $m$ satisfies 
\begin{equation}\label{eq 5}
\|T_{m}\|_{L^2\times L^2\times L^2\to L^{2/3}}\lesssim_{C_0,d,q}\|m\|^{q/3}_{L^{q}}.\end{equation}
Conversely, there is a function $m\in\bigcap_{q>3} L^{q}(\mathbb{R}^{3d})\cap\mathcal{L}^{\infty}(\mathbb{R}^{3d})$ such that the associate operator $T_{m}$ does not map $L^2\times L^2\times L^2\to L^{2/3}$. 

\section{Sufficiency}

\subsection{Decomposition scheme}

\noindent
The division scheme is as follows. Firstly, one starts with a fixed $j,G$ and defines $m^{j,G} := \sum_{n} b^{j,G}_{n}\Phi^{j,G}_{n}$ and further decomposes it to
\begin{equation*}m^{j,G} = \sum_{r}m^{j,G,r}.\end{equation*}
Altogether, 
\begin{equation*}m=\sum_{j,G}\sum_{r}m^{j,G,r}.\end{equation*}
One wishes to obtain through calculations that an expression of the following form,
\begin{equation*}\|T_{m}\|_{L^2\times L^2\times L^2\to L^{2/3}}\lesssim \sum_{j}\sum_{r}C(j,r)\prod_{i}\|f_{i}\|_{L^2}\end{equation*}
where $\sum_{j}\sum_{r}C(j,r)<\infty$. The summing over $j,r$ is over the expressed order: first $r$ then $j$. The sum in $G$ is implied in the summing over $j$, due to the definition of $\Phi^{j,G}_{n}$.

\subsection{Set-up}

\noindent
Let $j,G$ be as in $\ref{eq 1}$ and $n = (k_1,k_2,k_3)\in\mathbb{Z}^{3d}$. Set
\begin{equation*}b^{j,G}_{n} = \langle\Phi_{n}^{j,G},m\rangle.\end{equation*}
Then (see the appendix)
\begin{equation}\label{eq Y1}
\|m\|_{L^{q}}\asymp_{d}\bigg\|\bigg(\sum_{(j,G)}\sum_{n\in\mathbb{Z}^{3d}}|b^{j,G}_{n}2^{3jd/2}\chi_{Q_{jn}}|^2\bigg)^{1/2}\bigg\|_{L^{q}},\end{equation} 
with $Q_{jn}$ being a cube centered at $2^{-j}n$ with side-length $2^{1-j}$. Let $\tilde{Q}_{jn}=(1/2)Q_{jn}$. Temporarily fix $j,G$ and refer to $b^{j,G}_{n}$ as simply $b_{n}$. Then due to the pairwise disjoint of the cubes $\tilde{Q}_{jn}$ in $n$ (for fixed $j,G$), one also has from $\ref{eq Y1}$ that
\begin{align*}
\|m\|_{L^{q}} &\gtrsim 2^{3jd/2}\bigg\|\bigg(\sum_{n\in\mathbb{Z}^{3d}}|b_{n}|^2\chi_{Q_{jn}}\bigg)^{1/2}\bigg\|_{L^{q}}\\
&\geq 2^{3jd/2}\bigg\|\bigg(\sum_{n\in\mathbb{Z}^{3d}}|b_{n}|^2\chi_{\tilde{Q}_{jn}}\bigg)^{1/2}\bigg\|_{L^{q}}\\
&=2^{3jd/2}\bigg\|\sum_{n\in\mathbb{Z}^{3d}}|b_{n}|\chi_{\tilde{Q}_{jn}}\bigg\|_{L^{q}}\\
&= 2^{3jd(1/2-1/q)}\bigg(\sum_{n\in\mathbb{Z}^{3d}}|b_{n}|^{q}\bigg)^{1/q}.\end{align*}
The disjointness of the cubes $\tilde{Q}_{jn}$ is used in the last two equalities above. Let $b=(b_{n})_{n\in\mathbb{Z}^{3d}}$. Then the calculation above says, 
\begin{equation}\label{eq Y2}\|b\|_{l^{q}}\lesssim 2^{3jd(1/q-1/2)}\|m\|_{L^{q}}.\end{equation}

\noindent
Let $r\in\mathbb{N}_0$. Define,
\begin{equation*}
U_{r} =\{(k_1,\cdots,k_{l})\in\mathbb{Z}^{3d}: 2^{-r-1}\|b\|_{l^{\infty}}<|b_{(k_1,k_2,k_3)}|\leq 2^{-r}\|b\|_{l^{\infty}}\}.\end{equation*}
Let $\mathcal{C}$ denote the cardinality of $U_{r}$. Note that $\mathcal{C}$ is at most:
\begin{equation}\label{eq Cl1}
\mathcal{C}\asymp 2^{rq}\|b\big|_{U_{r}}\|_{l^{q}}^{q}\|b\|_{l^{\infty}}^{-q},\end{equation}
where $b\big|_{U_{r}}$ means $\{b_{n}\}_{n\in U_{r}}$. Define,
\begin{equation*}m^{r} := \sum_{(k_1,k_2,k_3)\in U_{r}}b_{(k_1,k_2,k_3)}\omega_{1,k_1}\omega_{2,k_2}\omega_{3,k_3} = \sum_{n\in U_{r}}b_{n}\omega_{1,k_1}\omega_{2,k_2}\omega_{3,k_3},\end{equation*} 
where the meaning of $\omega_{i,k_{i}}$ is as described in {\it Remark 4}, and correspondingly,
\begin{equation*}T_{m^{r}}(f_1,f_2,f_3)(x):= \int_{\mathbb{R}^{3d}} m^{r}(\xi_1,\xi_2,\xi_3)\hat{f}_1(\xi_1)\hat{f}_2(\xi_2)\hat{f}_3(\xi_3)e^{2\pi ix\cdot\sum_{i}\xi_{i}}\, d\xi_1d\xi_2 d\xi_3.\end{equation*}
In terms of Fourier transforms, $T_{m^{r}}$ can be written as,
\begin{equation*}T_{m^{r}}(f_1,f_2,f_3)(x)=\sum_{(k_1,k_2,k_3)\in U_{r}}b_{n}\mathcal{F}^{-1}(\omega_{1,k_1}\hat{f}_1)\mathcal{F}^{-1}(\omega_{2,k_2}\hat{f}_2)\mathcal{F}^{-1}(\omega_{3,k_3}\hat{f}_3).\end{equation*}

\noindent
The following lemma is key:\\

\noindent
{\bf Lemma 8.} For all $r\in\mathbb{N}$, let $\mathcal{C}$ be corresponding as above. Then
\begin{equation}\label{eq L1}
\|T_{m^{r}}(f_1,f_2,f_3)\|_{L^{2/3}}\lesssim_{d,q} 2^{3jd/2}(2^{-r}\|b\|_{l^{\infty}})\cdot\mathcal{C}^{1/3}\|f_1\|_{L^2}\|f_2\|_{L^2}\|f_3\|_{L^2}.\end{equation}

\subsection{Proof of Lemma 8} 

\noindent
Let
\begin{align*}
U_{r}^1 &=\{(k_1,k_2,k_3)\in U_{r}: card\{(k_1,k_2):n\in U_{r}\}>\mathcal{C}^{8/9}\},\\
U_{r}^2 &=\{(k_1,k_2,k_3)\in U_{r}\setminus U_{r}^1: card\{(k_2,k_3):n\in U_{r}\setminus U_{r}^1\}>\mathcal{C}^{8/9}\}.\\
U_{r}^3 &=\{(k_1,k_2,k_3)\in U_{r}\setminus(U_{r}^1\cup U_{r}^2): card\{(k_1,k_3):n\in U_{r}\setminus(U_{r}^1\cup U_{r}^2)\}>\mathcal{C}^{8/9}\},\end{align*}
and $U_{r}^4$ is the remaining set. 
Let $m^{r,i}$ denote the multipliers associated with $U^{i}_{r}$, respectively, and $T_{m^{r,i}}$ the operator associated with $m^{r,i}$, where
\begin{equation*}m^{r,i} =\sum_{n\in U^{i}_{r}}b_{n}\omega_{1,k_1}\omega_{2,k_2}\omega_{3,k_3}.\end{equation*}

\subsubsection{The operators $T_{m^{r,i}}$, $1\leq i\leq 3$}

\noindent
First consider $T^{m^{r,2}}$; others are similar. Let $E=\{k_1:(k_1,k_2,k_3)\in U^2_{r}\}$. \noindent
A simple calculation shows that
\begin{equation}\label{eq Q1}
|E|\leq\mathcal{C}^{1/9}.\end{equation}

\noindent
For each $k_1\in E$, denote
\begin{equation*}\mathcal{T}^{k_1,2}(f_2,f_3)=\sum_{(k_2,k_3): n\in U_{r}^2}b_{n}\prod_{i\geq 2}\mathcal{F}^{-1}(\omega_{i,k_{i}}\hat{f}_{i}),\end{equation*}
and,
\begin{equation*}A^{k_1}=\{(k_2,k_3): n=(k_1,k_2,k_3)\in U_{r}^2\}\end{equation*}
in the definition of $\mathcal{T}^{k_1,2}$. Let $\mathcal{A}^{k_1}=|A^{k_1}|$. Then the following is essentially done in [7]
\begin{equation}\label{eq Q2}
\|\mathcal{T}^{k_1,2}\|_{L^2\times L^2\to L^1}\lesssim_{d,q} 2^{3jd/2}(2^{-r}\|b\|_{l^{\infty}})\cdot\mathcal{A}_{k_1}^{1/4}.\end{equation}
The few things to take away here before applying the result in [7]. Firstly, the sequence $\{b_{n}\}_{n}$ in [7] is such that $n\in\mathbb{Z}^{2d}$. However in this case, with every fixed $k_1$, $b_{n}$ here can be regarded as having $n\in\mathbb{Z}^{2d}$, and the RHS dominance of $\ref{eq Q2}$ still holds. The conclusion in [7] is that the power of two should be $2^{jd}$ as it's for the bilinear case. But clearly $2^{jd}\leq 2^{3jd/2}$. Then from H\"older's inequality
\begin{align}\label{eq Q3}
\nonumber \|T_{m^{r,2}}(f_1,f_2,f_3)\|_{L^{2/3}}^{2/3} &\leq \int_{\mathbb{R}^{d}}\bigg(\sum_{k_1\in E} |\mathcal{F}^{-1}(\omega_{1,k_1}\hat{f}_1)|\cdot\bigg|\mathcal{T}^{k_1,2}(f_2,f_3)\bigg|\bigg)^{2/3}\,dx\\
\nonumber &\leq\bigg\|\bigg(\sum_{k_1\in E}|\mathcal{F}^{-1}(\omega_{1,k_1}\hat{f}_1)|^2 \bigg)^{1/2}\bigg\|_{L^2}^{2/3}\bigg\|\bigg(\sum_{k_1\in E}|\mathcal{T}^{k_1,2}(f_2,f_3)|^2\bigg)^{1/2}\bigg\|_{L^1}^{2/3}\\
&=: I_1\times I_2.\end{align}

\noindent
Now $I_1$ is simply,
\begin{multline*}
I_1= \bigg(\int_{\mathbb{R}^{d}}\sum_{k_1\in E}|\mathcal{F}^{-1}(\omega_{1,k_1}\hat{f}_1)|^2\,dx\bigg)^{1/3}\leq\bigg(\sum_{k_1\in E}\bigg(\int_{\mathbb{R}^{d}}|\mathcal{F}^{-1}(\omega_{1,k_1}\hat{f}_1)|^2\,dx\bigg)^{1/3}\\
=\bigg(\sum_{k_1\in E}\|\omega_{1,k_1}\hat{f}_1\|_{L^2}^2\bigg)^{1/3}\lesssim_{d} \bigg(\|\sum_{k_1\in E}\omega_{1,k_1}\hat{f}_1\|_{L^2}^2\bigg)^{1/3}\lesssim_{d}2^{jd/3}\|f_1\|_{L^2}^{2/3},\end{multline*}
from again, {\it Remark 4}, $\ref{eq 2}$. For $I_2$ in $\ref{eq Q3}$, one has, 
\begin{multline}\label{eq Q4} 
I_2= \bigg(\int_{\mathbb{R}^{d}}\bigg(\sum_{k_1\in E}|\mathcal{T}^{k_1,2}(f_2,f_3)|^2\bigg)^{1/2}\,dx\bigg)^{2/3}\\
\leq\bigg(\sum_{k_1\in E}\int_{\mathbb{R}^{d}}|\mathcal{T}^{k_1,2}(f_2,f_3)|\,dx\bigg)^{2/3} =\bigg(\sum_{k_1\in E}\|\mathcal{T}^{k_1,2}(f_2,f_3)\|_{L^1}\bigg)^{2/3}.\end{multline}
From $\ref{eq Q2}$ and the definition of $U_{r}^2$, the factor inside the root in $\ref{eq Q4}$ is dominated by,
\begin{equation*}\sum_{k_1\in E}\|\mathcal{T}^{k_1,2}(f_2,f_3)\|_{L^1} \lesssim_{d} \sum_{k_1\in E} 2^{jd}(2^{-r}\|b\|_{l^{\infty}})\mathcal{A}_{k_1}^{1/4}\prod_{i\geq 2}\|f_{i}\|_{L^2},\end{equation*}
where
\begin{align*}\sum_{k_1\in E}\mathcal{A}_{k_1}^{1/4} &=\bigg(\sum_{k_1\in E}\mathcal{A}_{k_1}^{1/2}\bigg)^{1/2}\bigg(\sum_{k_1\in E} 1\bigg)^{1/2}\\
& = \bigg(\mathcal{C}\mathcal{A}_{k_1}^{-1/2}\bigg)^{1/2}\bigg(\sum_{k_1\in E} 1\bigg)^{1/2}\\
&\leq (\mathcal{C}\mathcal{C}^{-4/9})^{1/2}(\mathcal{C}^{1/9})^{1/2} =\mathcal{C}^{1/3}.\end{align*}

\noindent
All of this implies,
\begin{equation}\label{eq Q5}
\|T_{m^{r,2}}(f_1,f_2,f_3)\|_{L^{2/3}}^{2/3}\lesssim_{d}2^{3jd/2}2^{-r}\|b\|_{l^{\infty}}\mathcal{C}^{1/3}\prod_{i}\|f_{i}\|_{L^2}\end{equation}
as needed.

\subsection{The operator $T_{m^{r,4}}$}

\noindent
The set $U_{r}^4$ has the following description:
\begin{equation*}U_{r}^4=\{(k_1,k_2,k_3):card\{(k_{i},k_{j}): n\in U_{r}\}\leq\mathcal{C}^{8/9}\}.\end{equation*}
Certainly $|U_{r}^4|\leq\mathcal{C}$. Divide $U_{r}^4$ into two sets, $V_1,V_2$ such that,
\begin{equation*}|V_1|,|V_2|\lesssim\mathcal{C}^{1/2}.\end{equation*}
One decomposes $T_{m^{r,4}}$ into two operators, each along different $V_{i}$. Consider $V_1$ first; the other is similar.\\

\noindent
Denote $V_1$ by $V$ and $|V|=\mathcal{V}$. Similarly as before, divide $V$ into four sets:
\begin{align*}
W_1 &=\{(k_1,k_2,k_3)\in V: card\{(k_1,k_2):n\in V\}>\mathcal{V}^{8/9}\},\\
W_2 &=\{(k_1,k_2,k_3)\in V\setminus W_1: card\{(k_2,k_3):n\in V\setminus W_2\}>\mathcal{V}^{8/9}\}.\\
W_3 &=\{(k_1,k_2,k_3)\in V\setminus (W_1\cup W_2): card\{(k_1,k_3):n\in V\setminus (W_1\cup W_2)\}>\mathcal{V}^{8/9}\},\end{align*}
and $W_4$ is the remaining set. Let $T_{W_{i}}$ be the corresponding operator to the set $W_{i}$. In other words, let $n\in W_{i}$ and $m_{W_{i}}=\sum_{n\in W_{i}}b_{n}\omega_{1,k_1}\omega_{2,k_2}\omega_{3,k_3}$ and
\begin{equation*}T_{W_{i}}(f_1,f_2,f_3)=\sum_{n\in W_{i}}b_{n}\prod_{j}\mathcal{F}^{-1}(\omega_{j,k_{j}}\hat{f}_{j}).\end{equation*}

\noindent
Follow the same strategy as before for $T_{W_{i}}, i=1,2,3$. One has that similar results as in $\ref{eq Q5}$
\begin{equation*}
\|T_{W_{i}}(f_1,f_2,f_3)\|_{L^{2/3}}^{2/3}\lesssim_{d}2^{3jd/2}2^{-r}\|b\|_{l^{\infty}}\mathcal{V}^{1/3}\prod_{i}\|f_{i}\|_{L^2}\leq2^{3jd/2}2^{-r}\|b\|_{l^{\infty}}\mathcal{C}^{1/3}\prod_{i}\|f_{i}\|_{L^2}.\end{equation*}

\noindent
For the last operator $T_{W_4}$ associated with the set $W_4$, one follows the same strategy of first a half dividing (by size) followed by a fourth dividing (by restriction on $2d$-slices) until one reaches the sets $S_{l}$ where each $S_{l}$ satisfies, 
\begin{equation*}|S_{l}|\leq\mathcal{C}^{1/8}.\end{equation*}
Pick one such set and call it $S$; the others are similar. The last step is to divide it into four sets, just as before,
\begin{align*}
R_1 &=\{(k_1,k_2,k_3)\in S: card\{(k_1,k_2):n\in S\}>|S|^{8/9}\},\\
R_2 &=\{(k_1,k_2,k_3)\in S\setminus R_1: card\{(k_2,k_3):n\in S\setminus R_1\}>|S|^{8/9}\}.\\
R_3 &=\{(k_1,k_2,k_3)\in S\setminus (R_1\cup R_2): card\{(k_1,k_3):n\in S\setminus (R_1\cup R_2)\}>|S|^{8/9}\}.\end{align*}
and $R_4$ is the remaining set. These sets by definitions, are disjoint. Let $T_{R_{i}}$ be the operator associated with $R_{i}$, just as before. Then the story is the same for $T_{R_{i}}, i=1,2,3$, and one has, 
\begin{equation*}
\|T_{R_{i}}(f_1,f_2,f_3)\|_{L^{2/3}}^{2/3}\lesssim_{d}2^{3jd/2}2^{-r}\|b\|_{l^{\infty}}\mathcal{C}^{1/3}\prod_{i}\|f_{i}\|_{L^2}.\end{equation*}
The last set $R_4$ has the following description:
\begin{equation*}R_4\subset\{(k_1,k_2,k_3)\in S\setminus (R_1\cup R_2\cup R_3): card\{(k_{i},k_{j}):n\in S\}\leq|S|^{8/9}\leq\mathcal{C}^{1/9}\}.\end{equation*}
That means, one can divide $R_4$ into at most $\mathcal{C}^{1/3}$ sets such that within each of these sets, if $(k_1,k_2,k_3)\not= (k_1',k_2',k_3')$ then $k_1\not=k_1', k_2\not=k_2',k_3\not= k_3'$. In other words, $(k_1,k_2,k_3)\leftrightarrow k_{i}$ is a one-to-one correspondence for each $i$. To see this, observe that for each fixed $k_1$, there exist at most $\mathcal{C}^{1/9}$ triples $(k_1,k_2,k_3),(k_1,k_2',k_3')$ in $R_4$, otherwise it would lead to a contradiction. This leads to at most $\mathcal{C}^{1/9}$ divisions in each of which $k_1$ components are all different. In each of these divisions, it might be the case that $(k_1,k_2,k_3)\not=(k_1',k_2',k_3')$ yet $k_2=k_2'$. Then further divides each of these divisions along the values of $k_2$'s, noting that by definition, it can be only at most $\mathcal{C}^{1/9}$ such divisions. The last division will be along the values of $k_3$. Altogether, there are at most $\mathcal{C}^{1/3}$ divisions or subsets, within each of which, the stated claim is satisfied. \\

\noindent
Let $D$ denote one of these sets and $T_{D}$ the operator associated with this set (defined similarly as before). Then by the generalized H\"older's inequality and the stated property of $D$,
\begin{align*}
\nonumber \|T_{D}(f_1,f_2,f_3)\|_{L^{2/3}}^{2/3} &\leq\sum_{n\in D}\|b_{n}\mathcal{F}^{-1}(\omega_{1,k_1}\hat{f}_1)\mathcal{F}^{-1}(\omega_{2,k_2}\hat{f}_2)\mathcal{F}^{-1}(\omega_{3,k_3}\hat{f}_3)\|_{L^{2/3}}^{2/3}\\
\nonumber &\leq (2^{-r}\|b\|_{l^{\infty}})^{2/3}\sum_{n\in D}\|\omega_{1,k_1}\hat{f}_1\|_{L^2}^{2/3}\|\omega_{2,k_2}\hat{f}_2\|_{L^2}^{2/3}\|\omega_{3,k_3}\hat{f}_3\|_{L^2}^{2/3}\\
\nonumber &\leq (2^{-r}\|b\|_{l^{\infty}})^{2/3}(\sum_{n\in D}\|\omega_{1,k_1}\hat{f}_1\|_{L^2}^2)^{1/3}(\sum_{n\in D}\|\omega_{1,k_1}\hat{f}_1\|_{L^2}^2)^{1/3}(\sum_{n\in D}\|\omega_{1,k_1}\hat{f}_1\|_{L^2}^2)^{1/3}\\
\nonumber &\leq (2^{-r}\|b\|_{l^{\infty}})^{2/3}(\|\sum_{k_1}\omega_{1,k_1}\hat{f}_1\|_{L^2}^2)^{1/3}(\|\sum_{k_2}\omega_{2,k_2}\hat{f}_2\|_{L^2}^2)^{1/3}(\|\sum_{k_3}\omega_{3,k_3}\hat{f}_3\|_{L^2}^2)^{1/3}\\
&\lesssim_{d} 2^{jd}(2^{-r}\|b\|_{l^{\infty}})^{2/3}\|f_1\|_{L^2}^{2/3}\|f_2\|_{L^2}^{2/3}\|f_3\|_{L^2}^{2/3}.\end{align*}

\noindent
Altogether,
\begin{equation*}\|T_{R_4}(f_1,f_2,f_3)\|_{2/3}\lesssim_{d}2^{3jd/3}(2^{-r}\|b\|_{l^{\infty}})\mathcal{C}^{1/3}\|f_1\|_{L^2}\|f_2\|_{L^2}\|f_3\|_{L^2}.\end{equation*}

\subsubsection{Conclusion}

\noindent
There are finite of division steps, whose number is independent of $\mathcal{C}$ (but depending on the trilinearity). Using the fact that 
\begin{equation*}\|f+g\|_{L^{2/3}}\leq C(\|f\|_{L^{2/3}}+\|g\|_{L^{2/3}}),\end{equation*}
for some $C>1$, over a finite number of times, one arrives at $\ref{eq L1}$,
\begin{equation*}\|T_{m^{r}}(f_1,f_2,f_3)\|_{2/3}\lesssim_{d}2^{3jd/3}(2^{-r}\|b\|_{l^{\infty}})\mathcal{C}^{1/3}\|f_1\|_{L^2}\|f_2\|_{L^2}\|f_3\|_{L^2}.\end{equation*}

\subsection{Conclusion for the sufficiency part}

\noindent
The {\bf Lemma 8} is helpful, because the next step in this analysis is to make sure that the obtained dominance for $\|T_{m^{r}}(f_1,f_2,f_3)\|_{L^1}$ in $\ref{eq L1}$ is summable over $r,j,G$ in that order, as explained before. Putting $\ref{eq Y2}, \ref{eq Cl1}$ into $\ref{eq L1}$ gives:
\begin{equation}\label{eq 22}
\|T_{m^{r}}\|_{L^2\times L_2\times L^2\to L^{2/3}}\lesssim_{d,q} 2^{jd(5/2-q/2)}2^{r(q/3-1)}\|b\|_{l^{\infty}}^{1-q/3}\|m\|_{L^{q}}^{q/3}.\end{equation}
Now $\ref{eq 22}$ is summable in $r$ as long as 
\begin{equation*}(q/3)<1\end{equation*}
which gives the range $1\leq q<3$. After summing $\ref{eq 22}$ in terms of $r$, the dominance constant on the RHS is reduced to:
\begin{equation*}
2^{jd(5/2-q/2)}\|b\|_{l^{\infty}}^{1-q/3}\|m\|_{L^{q}}^{q/3},\end{equation*}
which, after recalling from $\ref{eq 4}$ that $\|b\|_{l^{\infty}}\lesssim_{C_0} 2^{-(K+1+d)j}$, where $K+1$ is the presupposed smoothness degree of $m$, becomes,
\begin{equation}\label{eq 23}
2^{jd(5/2-q/2)}(2^{-(K+1+d)j})^{1-q/3}\|m\|_{L^{q}}^{q/3}.\end{equation}
Summing $\ref{eq 23}$ in terms of $j$ and requiring convergence gives,
\begin{equation*}(K+d+1)(1-q/3)>d\Rightarrow K+d+1>\frac{3d}{3-q}.\end{equation*}
That means one can require $M_{q} := K+1 = \lfloor\frac{3d}{3-q}\rfloor+1$, for instance.\\

\noindent
{\it Remark 6:} The dependence on the constant $C_0$ in $\ref{eq 5}$ is as follows. Note that $C_0$ was first invoked in {\bf Lemma 3}, which was done in [6]. In that paper, the dependence of $C_0$ is found to be linear. In particular, the RHS of $\ref{eq 4}$ can be taken to be $C_0 2^{-(K+1+d)j}$. This dominance is then put back in $\ref{eq 23}$, which then gives the total dominance for $\|T_{m}\|_{L^2\times L^2\times L^2\to L^{2/3}}$ as,
\begin{equation*}\|T_{m}\|_{L^2\times L^2\times L^2\to L^{2/3}}\lesssim_{d,q}C_0^{1-q/3}\|m\|_{L^{q}}^{q/3}.\end{equation*}

\section{Necessity}

\noindent
Let $\phi$ be a Schwartz function on $\mathbb{R}$ whose Fourier transform has support in a symmetric interval $I$ and let $\{a_{j}\}_{j\geq 1},\{b_{j}\}_{j\geq 1}, \{c_{j}\}_{j\geq 1}$ be two sequences of nonnegative numbers with only finitely many nonzero terms. This function $\phi$ is not related to the wavelet functions in the previous sections. Define $f,g,h$ by
\begin{equation}\label{eq 35}\hat{f}(\xi) = \sum_{j\geq 1}a_{j}\hat{\phi}(\xi_1-j)\prod_{r\geq 2}\hat{\phi}(\xi_{r}-1),\end{equation}
\begin{equation}\label{eq 36}\hat{g}(\eta) = \sum_{j\geq 1}b_{j}\hat{\phi}(\eta_1-j)\prod_{r\geq 2}\hat{\phi}(\eta_{r}-1),\end{equation}
\begin{equation}\label{eq 37}\hat{h}(\delta) = \sum_{j\geq 1}c_{j}\hat{\phi}(\delta_1-j)\prod_{r\geq 2}\hat{\phi}(\delta_{r}-1).\end{equation}
Then $f,g,h$ are Schwartz functions whose $L^2$ norms are bounded by a constant multiple of $(\sum_{j\geq 1}a_{j}^2)^{1/2}, (\sum_{j\geq 1}b_{j}^2)^{1/2}, (\sum_{j\geq 1}c_{j}^2)^{1/2}$, respectively.\\

\noindent
Let $\{s_{j}(t)\}_{j\geq 1}$ denote the sequence of Rademacher functions [3]. Let $\{v_{j}\}_{j\geq 1}$ be a bounded sequence of nonnegative numbers. For $t\in [0,1]$, consider $m_{t}$ by
\begin{equation}\label{eq 38}
m_{t}(\xi,\eta,\delta)=\sum_{j\geq 1}\sum_{k\geq 1}\sum_{l\geq 1}v_{j+k+l}s_{j+k+l}(t)\psi(\xi_1-j)\psi(\eta_1-k)\psi(\delta_1-l)\prod_{r\geq 2}\psi(\xi_{r}-1)\psi(\eta_1-1)\psi(\delta_1-1).\end{equation}
Here $\psi$ is a smooth function on $\mathbb{R}$ supported in the interval $J = 10 I$ and assuming value $1$ in $cJ$, with $c$ small enough so that $I\subset cJ$. Then from the definitions $\ref{eq 35}, \ref{eq 36}, \ref{eq 37}$,
\begin{multline}\label{eq 39} 
T_{m_{t}}(f,g,h)(x)\\
=\sum_{j\geq 1}\sum_{k\geq 1}\sum_{l\geq 1}a_{j}b_{k}c_{l}v_{j+k+l}s_{j+k+l}(t)\phi(x_1)^3e^{2\pi ix_1(j+k+l)}\prod_{r\geq 2}e^{6\pi ix_{r}}\phi(x_{r})^3\\
=\sum_{l\geq 3}v_{l}s_{l}(t)\phi(x_1)^3e^{2\pi ix_1l}\sum_{j=1}^{l-1}\sum_{k=1}^{l-j-1} a_{j}b_{k}c_{l-j-k}\prod_{r\geq 2}e^{6\pi ix_{r}}\phi(x_{r})^3.
\end{multline}

\noindent
Utilizing $\ref{eq 39}$, Khinchin's inequality [3] and Fubini's theorem, one then has
\begin{multline}\label{eq 40}
\int_0^1\|T_{m_{t}}(f,g,h)\|_{L^{2/3}}^{2/3}\,dt\\
=\bigg(\int_{\mathbb{R}}|\phi(y)|^2\,dy\bigg)^{n-1}\int_{\mathbb{R}}\int_0^1\bigg|\sum_{l\geq 3}v_{l}s_{l}(t)\phi(x_1)^3e^{2\pi ix_1l}\sum_{j= 1}^{l-1}\sum_{k=1}^{l-j-1} a_{j}b_{k}c_{l-j-k}\bigg|^{2/3}\,dt dx_1\\
\approx_{C,d}\int_{\mathbb{R}}\bigg(\sum_{l\geq 3}\bigg(v_{l}|\phi(x_1)|^3\sum_{j=1}^{l-1}\sum_{k=1}^{l-j-1} a_{j}b_{k}c_{l-j-k}\bigg)^2\bigg)^{1/3}dx_1 \approx_{C,d}\bigg(\sum_{l\geq 3}\bigg(v_{l}\sum_{j=1}^{l-1}\sum_{k=1}^{l-j-1} a_{j}b_{k}c_{l-j-k}\bigg)^2\bigg)^{1/3}.
\end{multline}

\noindent
Fix a positive integer $N\geq 2$ and set $a_{j}^{N}=b_{j}^{N}=c_{j}^{N}=2^{-N/2}$ if $j=2^{N},\cdots, 2^{N+1}-1$ and zero otherwise. Observe that with this agreement, $a_{j}b_{k}c_{l-j-k}$ is only nonzero for a finite number of terms for each $l$. Moreover $\sum_{j\geq 1}\sum_{k\geq 1}a_{j}b_{k}c_{l-j-k} =0$ if $l>2^{N+3}$ or if $l<2^{N+1}$. In effect, 
\begin{equation}\label{eq 41}\sum_{j=1}^{l-1}\sum_{k=1}^{l-j-1}a_{j}b_{k}c_{l-j-k}\geq\sum_{c2^{N}}^{C2^{N}}\sum_{c2^{N}}^{C2^{N}}2^{-3N/2}\geq c2^{N/2}>0.\end{equation}
In $\ref{eq 41}$ above, the constants $C,c$ in the first instance can be found as follows. Independently of $N$, $a_{j},b_{k}, c_{l-j-k}$ are "activated" (nonzero) for ranges of $j,k,l$ that have the same length. Hence within those ranges of $j,k$, one can fit in a "box" of sizes $[c2^{N},C2^{N}]\times [c2^{N},C2^{N}]$, for some appropriate $c,C$. The constant $c$ in the second instance of $\ref{eq 41}$ is unrelated to the previous one. \\

\noindent
Define $v_{l}=(l-1)^{-1}(\log(l-1))^{1/2}$ and define $f^{N}, g^{N}, h^{N}$ similarly to $f,g,h$, respectively, only with $a_{j},b_{j},c_{j}$ being replaced by $a_{j}^{N},b_{j}^{N},c_{j}^{N}$, respectively. Then from $\ref{eq 40}, \ref{eq 41}$ and simple calculus,
\begin{align*}
\int_0^1\|T_{m_{t}}(f,g,h)\|_{L^{2/3}}^{2/3}\,dt
&\gtrsim_{C,d} 2^{N/3}\bigg(\sum_{c2^{N}\leq l\leq C2^{N}}(l-1)^{-2}\log(l-1)\bigg)^{1/3}\\
&\gtrsim_{C,d} 2^{N/3}(\log c2^{N})^{1/3}\bigg(\sum_{c2^{N}\leq l\leq C2^{N}}(l-1)^{-2}\bigg)^{1/3}\\
&\gtrsim 2^{N/3}cN^{1/3}\bigg(\int_{c2^{N}}^{C2^{N}} y^{-2}\,dy\bigg)^{1/3}\gtrsim cN^{1/3}.\end{align*}
Above, $c,C$ denote positive numbers that possibly change from one instance to the next, with $C>c$. This means that for every $N\geq 2$ one can find $t_{N}\in [0,1]$ such that
\begin{equation}\label{eq 42}
\|T_{m_{t_{N}}}(f,g,h)\|_{L^{2/3}}\geq CN^{1/2}\end{equation}
with $C$ being independent of $N$.\\

\noindent
Define:
\begin{equation*}
m(\xi,\eta,\delta)=\sum_{j\geq 1}\sum_{k\geq 1}\sum_{l\geq 1}v_{j+k+l}\sigma_{j+k+l}\psi(\xi_1-j)\psi(\eta_1-k)\psi(\delta_1-l)\prod_{r\geq 2}\psi(\xi_{r}-1)\psi(\eta_1-1)\psi(\delta_1-1)\end{equation*}
with $\sigma_{l} = s_{l}(t_{N})$ if $N\geq 2$ and $2^{N+1}\leq l\leq 2^{N+3}$. Then a quick calculation shows that,
\begin{equation*}
T_{m}(f^{N},g^{N},h^{N})(x)=T_{m_{t_{N}}}(f^{N},g^{N},h^{N})(x),\end{equation*}
which means, from $\ref{eq 42}$,
\begin{equation*}
\|T_{m}(f^{N},g^{N},h^{N})\|_{L^{2/3}}\gtrsim N^{1/2}.\end{equation*}
Hence $T_{m}$ is unbounded on $L^2\times L^2\times L^2\to L^{2/3}$. On the other hand with the definition, $m$ is a smooth function with bounded derivatives and moreover,
\begin{equation*}
\|m\|_{L^{q}}\approx_{q}\bigg(\sum_{l\geq 3} v_{l}^{q}(l-1)^2\bigg)^{1/q}=\bigg(\sum_{l\geq 3}(\log(l-1))^{q/2}(l-1)^{2-q}\bigg)^{1/q}\end{equation*}
which, by comparison with the divergent series $\sum_{n\geq 3} (1/n)$, shows that $m\in\bigcap_{q>3} L^{q}$ and $m\not\in L^{q}$ if $q\leq 3$.

\section{Appendix}

\noindent
The following utilizes the material in [9]. First one should see the definition of the function spaces $F^{s}_{pq}(\mathbb{R}^{d})$, $s\in\mathbb{R}$, $0<p\leq\infty, 0<q\leq\infty$ [9]. The important point here is,
\begin{equation} \label{eq 24} F^0_{q2} = L^{q}.\end{equation}
See the remark 1.65 in [9]. One also has the following fact, which is a rephrasing of parts of Theorem 1.64 in [9], using the notations already introduced in this paper:\\

\noindent
{\bf Theorem A.} Let $\Phi^{j,G}_{n}$ be as in $\ref{eq 1}$ with sufficient smoothness (see {\it Remark 3}). Let $f\in F^{s}_{qp}(\mathbb{R}^{d})$. Then $f$ can be represented as,
\begin{equation*} f= \sum_{(j,G,n)}b^{j,G}_{n}\Phi^{j,G}_{n}.\end{equation*}
Furthermore, this representation is unique and the map
\begin{equation*}I: f\mapsto\{b^{j,G}_{n}\}\end{equation*}
is an isomorphic map of $F^{s}_{qp}$ onto $f^{s}_{qp}$. The latter is a sequence space whose elements $\lambda = \{\lambda\}^{j,G}_{n}$ are given the following semi-norm:
\begin{equation}\label{eq 25}
\|\lambda\| = \bigg\|\bigg(\sum_{(j,G,n)}|2^{js}2^{jn}\lambda^{j,G}_{n}\chi_{Q_{jn}}|^{p}\bigg)^{1/q}\bigg\|_{L^{q}(\mathbb{R}^{d})}.\end{equation}
{\bf Theorem A}, $\ref{eq 24},\ref{eq 25}$ imply $\ref{eq Y1}$.

\subsection{Lemma 3}

\noindent
It was mentioned in [7] that the proof of {\bf Lemma 3} is essentially done in Lemma 7 of [6]. Since the proof in [6] is quite involved, only its ideas will be presented here. This discussion aims to guarantee that although the version of the lemma stated in [7] was stated for $2d$, it's immaterial to change it to any $ld$. \\

\noindent
First, a small point here is that, if one follows the argument there then because of the change in dimensions one should arrive at
\begin{equation*} |\langle\Phi^{j,G}_{n},m\rangle|\lesssim_{C_0} 2^{-(K+1+3d/2)j}\end{equation*}
for the conclusion. But clearly, $2^{-3d/2}\leq 2^{-d}$.\\

\noindent
If the multiplier had any decay of the form,
\begin{equation}\label{eq 26}
|\partial^{\alpha}m(\vec{x})|\leq (1+|\vec{x}|)^{M_1}\end{equation}
for $|\alpha|\leq K$, for some large $K, M_1$, then one can utilize this and the moment cancellation properties $\Phi^{j,G}_{n}$ ($j\not=0$) and apply the material in the Appendix B.2 of [4] to get the desired decay $\ref{eq 4}$ of the wavelet transforms of $m$. The moment cancellation properties of $\Phi^{j,G}_{n}$ when $j\not=0$ come from those of $\phi_{M}$ in its definition. There are no cancellation assumptions when $j=0$, but then one can use the decay $\ref{eq 24}$ and those of the wavelets to apply the material in Appendix B.1 of [4] instead. \\

\noindent
Our $m$ is assumed of no decay $\ref{eq 26}$. Hence one can decompose $m$ into parts, 
\begin{equation*}m =\sum_{i} m_{i}\end{equation*} 
with $m_{i}$ being defined as in [6]. Then $m_{i} = m_0(2^{i}\cdot)$. Then $m_0$ does possesses the decay $\ref{eq 26}$ [6]. One can arrive at $\ref{eq 4}$ for $m_0$. Now for $i<0$, one has through change of variables,
\begin{equation}\label{eq 27}
\langle\Phi^{j,G}_{n},m_0(2^{i}\cdot)\rangle = 2^{ic}\langle\Phi^{j+i,G}_{n},m_0\rangle =2^{-c|i|}\langle\Phi^{j+i,G}_{n},m_0\rangle .\end{equation}
Here $c=c(d) = d/2$ for any dimension $d$ of $\vec{x}$. Then one can take $j\geq 0$ in $\ref{eq 27}$ large enough so that $k=j+i\in\mathbb{N}_0$. For $i>0$, one uses, 
\begin{equation*}
\langle\Phi^{j,G}_{n},m_0(2^{i}\cdot)\rangle = 2^{-c|i|}\langle\Phi^{j-i,G}_{n},m_0\rangle.\end{equation*}
Putting all these back, one gets the desired conclusion for $m$. None of the tools mentioned here is particular to any dimension $d$. Hence the conclusion $\ref{eq 4}$ holds for $3d$ or any $ld$.

%\clearpage

\end{document}